\documentclass[letterpaper, 10 pt, conference]{ieeeconf}  

\IEEEoverridecommandlockouts                              

\usepackage{graphicx}
\usepackage{amsmath} 
\usepackage{amssymb}  
\usepackage{color}
\usepackage{cite}
\DeclareMathOperator{\e}{\mathrm{e}}
\newtheorem{definition}{Definition}
\newtheorem{proposition}{Proposition}
\newtheorem{lemma}{Lemma}
\newtheorem{theorem}{Theorem}
\newtheorem{corollary}{Corollary}

\title{\LARGE \bf
Stability Properties of the Impulsive Goodwin's Oscillator in 1-cycle
}

\author{Anton V. Proskurnikov$^{1}$ and Alexander Medvedev$^{2}$
\thanks{Anton V. Proskurnikov  [{\tt\small anton.p.1982@ieee.org}] is with the Department of Electronics and Telecommunications, Politecnico di Torino, Turin, Italy, 10129.}%
\thanks{Alexander Medvedev [{\tt\small alexander.medvedev@it.uu.se}] is with Department of Information Technology,
        Uppsala University, SE-752 37 Uppsala, Sweden.
        }%
}

\begin{document}

\maketitle
\thispagestyle{empty}
\pagestyle{empty}

\begin{abstract}
The Impulsive Goodwin's Oscillator (IGO) is a mathematical model that represents a hybrid closed-loop system. It arises by closing a specific of type continuous positive linear time-invariant system with impulsive feedback, incorporating both amplitude and frequency pulse modulation. The structure of the IGO precludes the existence of equilibria and thus ensures that all of its solutions, whether periodic or non-periodic, are oscillatory. Originating in mathematical biology, the IGO also constitutes a control paradigm applicable to a wide range of fields, particularly to closed-loop dosing of chemicals and medicines. The pulse modulated feedback  introduces strong nonlinearity and non-smoothness into the closed-loop dynamics thus rendering conventional controller design methods  not applicable. However, the hybrid dynamics of IGO reduce to a nonlinear discrete-time system, exhibiting a one-to-one correspondence between  solutions of the original hybrid IGO and those of the discrete-time system. The paper proposes a design approach that leverages the linearization of the equivalent discrete-time  dynamics in the vicinity of a fixed point. An efficient local stability condition of the 1-cycle in terms of the characteristics of the amplitude and frequency modulation functions is obtained.
Unlike the conventional Schur-Cohn and Jury stability conditions applied to the Jacobian matrix, the obtained criterion requires checking a single inequality that is linear in the slopes of the modulation characteristics.
\end{abstract}

\section{INTRODUCTION}

Most research in control theory focuses on steering a dynamical system towards and stabilizing it at an equilibrium point. However, there is growing interest in oscillatory behaviors that are ubiquitous  in physics, chemistry, biology, economics, engineering, and medicine~\cite{LRS22}. Modeling and analysis of periodic and non-periodic oscillations are therefore timely topics in nonlinear dynamics, with rich applications across science and technology.

A periodic oscillation describes a process that repeats in a cyclic manner. A standard example of control  actions performed repeatedly according to a schedule is taking prescription drugs. Under stationary conditions, administering the right dose at the appropriate time usually works well. When the therapeutic effect is insufficient, either the drug dose must be increased or the interdose interval must be reduced. Both regimen adjustments elevate the drug concentration in the organism and eliciting a stronger effect according to the dose-response relationship.  Similar control mechanisms that adjust both the timing and magnitude of discrete actions also appear outside medicine, e.g., in mechanical systems with impacts~\cite{BB00} or in pest management~\cite{SC05}.

A control action that is fast relative to the plant dynamics can be approximated by an impulse, or the Dirac $\delta$-function. A well-developed framework to handle a continuous system with impulsive output feedback is pulse-modulated control with amplitude and frequency modulation~\cite{GC98,SP}. This framework was successfully utilized in the Impulsive Goodwin's Oscillator (IGO), devised to model pulsatile endocrine regulation \cite{MCS06,Aut09}. The IGO has found application in modeling biological data pertaining to feedback (non-basal) testosterone regulation in males~\cite{PMAM14}, pulsatile secretion of cortisol \cite{Runvik:2022}, and the multi-peak phenomenon in \emph{levodopa}, a drug used to treat Parkinson's disease~\cite{Runvik:2020}.

Since the inception of the IGO, the research has primarily focused on exploring the complex dynamical phenomena exhibited by the model, namely periodic solutions of high multiplicity, chaos \cite{ZCM12b}, and the entrainment of oscillations to an exogenous  periodic signal~\cite{MPZh18}. These studies pertain to the \emph{analysis} of the IGO's complex dynamics.

More recently, the \emph{design} of the pulse-modulated feedback of the IGO to sustain a desired periodic solution has been addressed. Two problems were solved for a given continuous plant, both concerned with a so-called 1-cycle -- a periodic solution characterized by a single firing of the feedback within the least period. First, in \cite{MPZh23}, the problem of obtaining a stable 1-cycle with a specific period and weight of the impulsive control sequence is solved. Second, the output corridor control of the continuous plant, which keeps the output within a predefined closed interval of values, was developed in \cite{MPZ24}. In both cases, applications to the dosing of chemicals and drugs were envisioned, e.g., dosing of a neuro-muscular blockade agent~\cite{MPZ24_ECC,MPZ24_MED} in general anesthesia.


This paper addresses the stability analysis of the 1-cycle in the IGO, which is essential for achieving a sustained periodic solution and is important in both analysis and design.
A simple and efficient local stability condition for the 1-cycle is derived. While this criterion is equivalent to the Schur stability of a Jacobian matrix, it requires validating only a single inequality that is linear in the slopes of the frequency and amplitude modulation characteristics. In solving the design problem, this condition provides an exact characterization of the IGO modulation functions   that ensure orbital stability of the desired closed-loop solution. Although this condition implies only local stability, extensive  numerical experiments show that, in fact, almost all solutions of the original hybrid IGO system are attracted to a 1-cycle if it is (locally) orbitally stable~\cite{ZCM12b}. Multistability has been observed only in the IGO with time delay~\cite{ZMCM15} and in the IGO driven by a continuous exogenous signal~\cite{MPZ18}.

The rest of the paper is composed as follows.
In Section~\ref{sec.igo}, the IGO model is introduced  for the reader's convenience. In this section, the dynamics of the IGO are discussed and a discrete map propagating the state vector of the continuous part of the model through the firings of the pulse-modulated feedback is given. The latter is used to derive an explicit expression for the 1-cycle. In Section~\ref{sec:stability}, the main result of the paper is formulated yielding a linear inequality that provides a necessary and sufficient condition of a 1-cycle with given parameters. Finally, an involved numerical example is considered in Section~\ref{sec:example} to illustrate how the developed theory informs the pulse-modulated controller design.

\section{IMPULSIVE GOODWIN'S OSCILLATOR}\label{sec.igo}
Consider a third-order linear time-invariant system
\begin{equation}                            \label{eq:1}
\dot{x}(t) =Ax(t), \quad y(t)=Cx(t).
\end{equation}
Here matrices $A,C$ are as follows:
\begin{equation}                            \label{eq:1a}
A=\begin{bmatrix} -a_1 &0 &0 \\ g_1 & -a_2 &0 \\ 0 &g_2 &-a_3 \end{bmatrix}, 
C =\begin{bmatrix}0 & 0 & 1\end{bmatrix},
\end{equation}
where $a_1,a_2,a_3>0$ are \emph{distinct} constants, and $g_1,g_2>0$ are positive gains. The scalar function $y$ is the measured output, and the state variables are
$x= [x_1,x_2,x_3]^\top$. It follows  that the matrix $A$ is Hurwitz stable and Metzler.

\subsection*{Impulsive feedback}

Continuous-time system~\eqref{eq:1} is controlled by a pulse-modulated feedback whith the impulse weights and their timing  determined by the continuous plant output $y(t)$:
\begin{align}\label{eq:2}
x(t_n^+) &= x(t_n^-) +\lambda_n B, \quad
t_{n+1} =t_n+T_n,
\\
T_n &=\Phi(y(t_n)), \; \lambda_n= F(y(t_n)), \quad B=\begin{bmatrix} 1 & 0 & 0\end{bmatrix}^{\top},\notag
\end{align}
where  $n=0,1,\ldots$.
The minus and plus in a superscript in~\eqref{eq:2} denote the left-sided and
a right-sided limit, respectively. Notice that the jumps in $x(t)$ lead  to discontinuities only in $x_1(t)$, whereas $x_2(t)$, and $y(t)=x_3(t)$ remain continuous.
The instants $t_n$ are called (impulse) firing times
and $\lambda_n$ represents the corresponding impulse weight.

In theory of pulse-modulated systems~\cite{GC98}, $F(\cdot)$ is called the amplitude modulation function and $\Phi(\cdot)$ is referred to as the frequency modulation function. The modulation functions are assumed to be  continuous and monotonic, $F(\cdot)$ be \emph{non-increasing}, and $\Phi(\cdot)$ be \emph{non-decreasing}.

These monotonicity assumptions imply that  controller~\eqref{eq:2} implements a negative feedback from the continuous output to the amplitude and frequency of the pulses. Namely,
an increased value of $y(t_n)$ results in a lower (or unchanged) weight $\lambda_n$ for the next impulse fired at $t_{n+1}$. Furthermore, the interval between the impulses increases, which makes the sequence of pulses sparser. This feedback mechanism prevents the controlled output from diverging, being similar to the continuous Goodwin's oscillator, where the feedback from the output $y(t)$ is implemented by means of a decreasing static nonlinearity~\cite{Good65}.
Notably, the negative feedback action is implemented by means of  positive signals only.

To explicitly restrict the domain where the solutions of closed-loop system \eqref{eq:1}, \eqref{eq:2} ultimately evolve,  boundedness of the modulation functions is
required
\begin{equation}                             \label{eq:2a}
\Phi_1\le \Phi(\cdot)\le\Phi_2, \quad 0<F_1\le F(\cdot)\le F_2,
\end{equation}
where $\Phi_1$, $\Phi_2$, $F_1$, $F_2$ are positive constants. Under these limitations, all solutions of \eqref{eq:1}, \eqref{eq:2} are contained in an invariant 3-dimensional box that can be computed explicitly~\cite{ZCM12b}.

\begin{definition}
The IGO is a hybrid system arising as a feedback interconnection of the continuous linear time-invariant block in \eqref{eq:1} and the impulsive feedback possessing first-order discrete dynamics in \eqref{eq:2}.
\end{definition}

The class of design problems that captures our interest involves guaranteeing certain desired properties of  the solutions to the IGO through the selection of the modulation functions $F(\cdot)$, $\Phi(\cdot)$. These functions serve as the designer's degrees of freedom in the impulsive controller tuning.

\subsection*{The discrete-time representation and 1-cycles}

The hybrid dynamics of the IGO can be reduced to a discrete-time equation by noticing that the sequence of state vectors
$X_n=x(t_n^-)$ obeys the recurrence formula~\cite{Aut09}
\begin{align}\label{eq:map}
    X_{n+1}&=Q(X_n),\\
    Q(\xi) &\triangleq \mathrm{e}^{A\Phi(C\xi)}\left( \xi+ F(C\xi)B \right).\nonumber
\end{align}
Since the plant is autonomous in between the impulsive feedback firings, the continuous state trajectory on the interval $(t_n,t_{n+1})$ is uniquely defined by $X_n$ as
\begin{equation} \label{eq:1d}
x(t)=\e^{(t-t_n)A}(X_n+\lambda_n B),\quad t\in(t_n,t_{n+1}).
\end{equation}
In this sense, the properties of the IGO, as a hybrid dynamical system, are completely determined by the properties of the impulse-to-impulse map $Q$, defined in~\eqref{eq:map}.

As has been reported in~\cite{ZCM12b}, discrete-time system~\eqref{eq:map} and, therefore, the IGO, can exhibit a wide range of periodic and non-periodic oscillation, including cycles of high multiplicity and deterministic chaos. In this study, only the simplest periodic solution of \eqref{eq:1}, \eqref{eq:2}  with one firing of the impulsive feedback in the least solution period is treated. It is termed 1-cycle, see e.g. \cite{MZ03}, and, by definition,
corresponds to the \emph{periodic} instants of pulses $t_{n+1}=t_n+T$, $T>0$ and the \emph{constant} sequence of amplitudes $X_n=X(t_n^-)\equiv X$, where $X$ is the fixed point of map $Q$:
\begin{equation}\label{eq:1-cycle}
    X=Q(X),
\end{equation}
The characteristics of the 1-cycle, i.e. the (least) period and the impulse weight,  are then defined by the fixed point as $T=\Phi(y_0)$, $\lambda=F(y_0)$, $y_0=CX$.

As previously demonstrated in~\cite{MPZh23}, the solution to  nonlinear equations~\eqref{eq:1-cycle} can be analytically expressed using the parameters of plant~\eqref{eq:1} and the characteristics of the 1-cycle $\lambda,T$. This explicit solution of~\eqref{eq:1-cycle} is conveniently formulated in terms of  divided differences (DDs). The first DD of a function $h:\mathbb{R}\to\mathbb{R}$ is defined as
\[
h\lbrack x_1, x_2 \rbrack \triangleq \frac{h(x_1)-h(x_2)}{x_1-x_2},\quad\forall x_1\ne x_2,
\]
and higher-order divided DDs are introduced recursively by
\[
h\lbrack x_0, \dots, x_k \rbrack = \frac{h\lbrack x_1, \dots, x_k\rbrack-h\lbrack x_0, \dots, x_{k-1}\rbrack }{x_k-x_0}.
\]
For the sake of simplicity,  only pairwise distinct sets of values $x_0,\ldots,x_k$ are considered here.

Denote $\mu(x)\triangleq\frac{1}{\e^{-x}-1}$. Using the standard definition of an analytic function on matrices~\cite{E87}, $\mu(M)=(\e^{-M}-I)^{-1}$ for an arbitrary non-singular matrix $M$.
The special structure of matrix $A$ defined in \eqref{eq:1a} allows to compute $\mu(TA)$, $T>0$, through DDs of the function $\mu$ by means of the Opitz formula~\cite{E87,PRM24} and
leads to the following proposition.

\begin{proposition}{\cite{MPZh23}}\label{pro:fp}
    If  IGO~\eqref{eq:1}, \eqref{eq:2} exhibits a 1-cycle of the period $T$ with the weight $\lambda$, then the fixed point  satisfying \eqref{eq:1-cycle} is uniquely determined as
    \begin{equation}\label{eq:fp_alpha}
        X= \lambda\mu(TA)B=\lambda  \begin{bmatrix}
            \mu(- a_1T) \\ g_1\mu\lbrack - a_1T,- a_2 T\rbrack \\  g_1g_2\mu\lbrack -a_1T, - a_2T, - a_3T \rbrack
        \end{bmatrix}.
    \end{equation}
\end{proposition}

The availability of an analytic expression for $X$ ensures one-to-one map between the pair $(\lambda,T)$ and the fixed point. This fact has enabled the \emph{design} of the IGO whose 1-cycles have predefined parameters $\lambda,T$~\cite{MPZh23,MPZh23a,MPZ24}. In order to guarantee the existence of such a 1-cycle, one has to find the (nonlinear) modulation functions $\Phi,F$ such that
\[
\begin{gathered}
\lambda=F(y_0),\quad T=\Phi(y_0),\;\;\text{where}\\
y_0\triangleq CX=\lambda g_1g_2\mu\lbrack -a_1T, - a_2T, - a_3T \rbrack.
\end{gathered}
\]
The problem is, however, that the resulting 1-cycle can turn out (orbitally) \emph{unstable} and thus fail to pertain in the face of perturbation. This is due to the fact that orbital stability of solutions to the IGO depends on the slopes of the modulation functions and not captured by the equations above.

In the next section, the main result of this paper is presented, offering a simple analytic stability criterion.

\section{STABILITY OF 1-CYCLE}\label{sec:stability}

The  1-cycle in closed-loop system \eqref{eq:1}, \eqref{eq:2}  corresponding to the fixed point $X$ is known to be (locally exponentially) orbitally stable \cite{Aut09,ZCM12b} if only only if
\begin{gather}\label{eq:jacobian}
    Q^\prime(X)= \mathrm{e}^{A\Phi(y_0)}\left( I+ F^\prime(y_0)BC\right)+ \Phi^\prime(y_0)AX C,\notag
\end{gather}
is a Schur stable matrix\footnote{Recall that a matrix is Schur stable, or Schur, if all its eigenvalues are less than $1$ in modulus.}. As pointed out in \cite[Proposition~3]{MPZh23},  the Jacobian can be written as
\begin{equation}\label{eq:closed_loop}
    Q^\prime(X)= \e^{A\Phi(y_0)}+ \begin{bmatrix}
    J &D
\end{bmatrix}\begin{bmatrix}
 F^\prime(y_0)  \\  \Phi^\prime(y_0)
\end{bmatrix}C,
\end{equation}
where $J=\e^{AT}B>0, D=AX<0$.

Since plant~\eqref{eq:1} is Hurwitz, stability of the 1-cycle
is always guaranteed for zero slopes of the modulation functions, for instance, when $F(y)=\mathrm{const}$, $\Phi(y)=\mathrm{const}$. However, this essentially eliminates the output feedback, at least in the vicinity of the fixed point. To improve the convergence to the stationary solution under perturbation, the spectral radius of the Jacobian has to be minimized.

\subsection*{Insufficiency of standard stabilization methods}

The right-hand side of \eqref{eq:closed_loop} has apparent similarity to the problem of stabilization of a discrete time-invariant linear system by a static output feedback, see e.g. \cite{CLS98}.  This is the problem of finding such gain matrix $K_d$ that the system
\begin{align*}
    x_d(t+1)&=A_d x_d(t)+ B_d u_d(t),\\
    y_d(t)&= C_dx_d(t),
\end{align*}
is (asymptotically) stabilized by the control law $ u_d(t)= K_d y_d(t)$. Equivalently, one is looking for a $K_d$ that makes the state matrix of the closed-loop system $A_d+B_dK_dC_d$ Schur-stable.
For the reasons described above, the largest possible set of such controllers is sought, despite the existence of the trivial solution $K_d\equiv 0$.

Although the static output feedback stabilization problem appears to be simple, a complete characterization of the gains solving it is missing. For instance, the pole placement problems via static feedback are usually considered in the situation where the total number of scalar entries $\dim y_d\dim u_d$ in $K_d$ is not less than the state dimension $\dim x_d$~\cite{Brockett81}, which inequality is, obviously, violated in the present case (cf. $\dim x_d=3,\dim u_d=2,\dim y_d=1$). This fact relates the problem at hand with the stabilization of underactuated systems that commonly appear in mechatronics and robotics \cite{ABF09}, where the control problem is solved via generating a sustained oscillation of a desired amplitude and frequency.

Another idea suggested by the similarity between the problem of stabilizing the fixed point of a 1-cycle  in \eqref{eq:1}, \eqref{eq:2}  and solving the static output feedback problem is to reformulate the stability condition as a system of bilinear matrix inequalities (BMI)
\begin{equation}\label{eq:BMI}
    (A_\Phi+WKE)^\top P (A_\Phi+WKE)-P<0,\;\;P>0,
\end{equation}
where
\[
A_\Phi=\e^{A\Phi(y_0)}, W=\begin{bmatrix}
    J &D
\end{bmatrix}, K= \begin{bmatrix}
  F^\prime(y_0)   & \Phi^\prime(y_0)
\end{bmatrix}^{\top},
\]
and $P$, $K$ are the decision variables. Again, one may notice that the inequality is feasible since  it is always satisfied for $K=0$ and some $P$.
However, the non-convexity of~\eqref{eq:BMI} makes it difficult to find the optimal (with respect to some performance index) solution; all known methods, in general, return only local solutions~\cite{Hassibi99}.

\subsection*{Main result: a linear stability condition}

Although the expression for the Jacobian matrix  is complicated, a necessary and sufficient analytic criteria for $Q'(X)$ being Schur stable is obtained below in terms of a \emph{linear inequality} in the slopes $F'(y_0),\Phi'(y_0)$. This makes the result here quite different from the standard Schur stability criteria, such as the Schur-Cohn stability test, the Jury criterion, and the Liénard-Chipart criterion (all of them lead lead to nonlinear stability conditions, see, e.g.,~\cite[Lemma~1]{MPZh23}). The next theorem is our main result, proved in Appendix.

\begin{theorem}\label{thm.main}
Assume that $0<a_1<a_2<a_3$.
If  $F^\prime(z_0)\leq 0$ and $\Phi^\prime(z_0)\geq 0$, then the Jacobian matrix $Q^\prime(X)$ is Schur stable if and only if
\begin{equation}\label{eq.stability-det}
\det(-I-Q^\prime(X))<0,
\end{equation}
or, equivalently, the following linear inequality holds
\begin{equation}\label{eq.stability-slopes}
C(I+e^{\Phi(y_0)A})^{-1}\left(F'(y_0)J+\Phi'(y_0)D\right)>-1.
\end{equation}
Furthermore, the Jacobian matrix $Q'(X)$ always has a positive real eigenvalue, lying in the interval $[\e^{-a_3T},\e^{-a_1T}]$. Hence, the spectral radius of $Q'(X)$ is not less than
$\e^{-a_3T}$.
\end{theorem}

\subsection*{Discussion}

Several insights can be gained from the analysis above. First,  since  $JF^\prime(y_0) + D\Phi^\prime(y_0)\leq 0$, for all feasible values of $F^\prime(y_0),\Phi^\prime(y_0)$, the feedback stabilizing the fixed point in \eqref{eq:closed_loop} is negative, in a well-defined sense, cf. \eqref{eq:closed_loop}. This is despite the fact that all the signals comprising closed-loop system \eqref{eq:1}, \eqref{eq:2} are positive.  Thus, the IGO mathematically explains how negative feedback is implemented in nature by impulsive regulation when negative signals are not available. Impulse-modulated feedback is  particularly important in the regulation of endocrine systems, e.g. in the hypothalamic-pituitary adrenal and gonadal axes, where amplitude and frequency modulation permits the encoding of signaling information to the target sites \cite{WTT10}.

Second,  existing design methods for (discrete) linear time-invariant systems can be adopted to  the framework of  IGO design in 1-cycle by  applying them to the fixed point instead of the equilibrium. This intriguing opening is unexpected given the complex nonlinear dynamics of the IGO. Yet, as pointed out in \cite{ZCM12b} with respect to the modeling of testosterone regulation in the male, the IGO with biologically reasonable parameter values exhibits mostly 1- or 2-cycles, which is in line with biological data \cite{PMAM14}.

Third, as a control structure, the IGO has a limit to the achievable local convergence rate due to the lower bound on the spectral radius imposed by the continuous plant. This is expected from the very principle of operation of the IGO as the dynamics of the closed-loop system are governed by \eqref{eq:1} in between the feedback firings.

Finally, the terms in \eqref{eq.stability-slopes} characterize the respective contributions of the amplitude and frequency modulation to achieving the control goal. This kind of characterization is instructive in e.g. the optimization of drug dosing regimen to decide on the balance between the dose and inter-dose time interval adjustments.

\section{NUMERICAL EXAMPLE}\label{sec:example}
\begin{figure}[h!]
\centering
\includegraphics[width=0.95\linewidth]{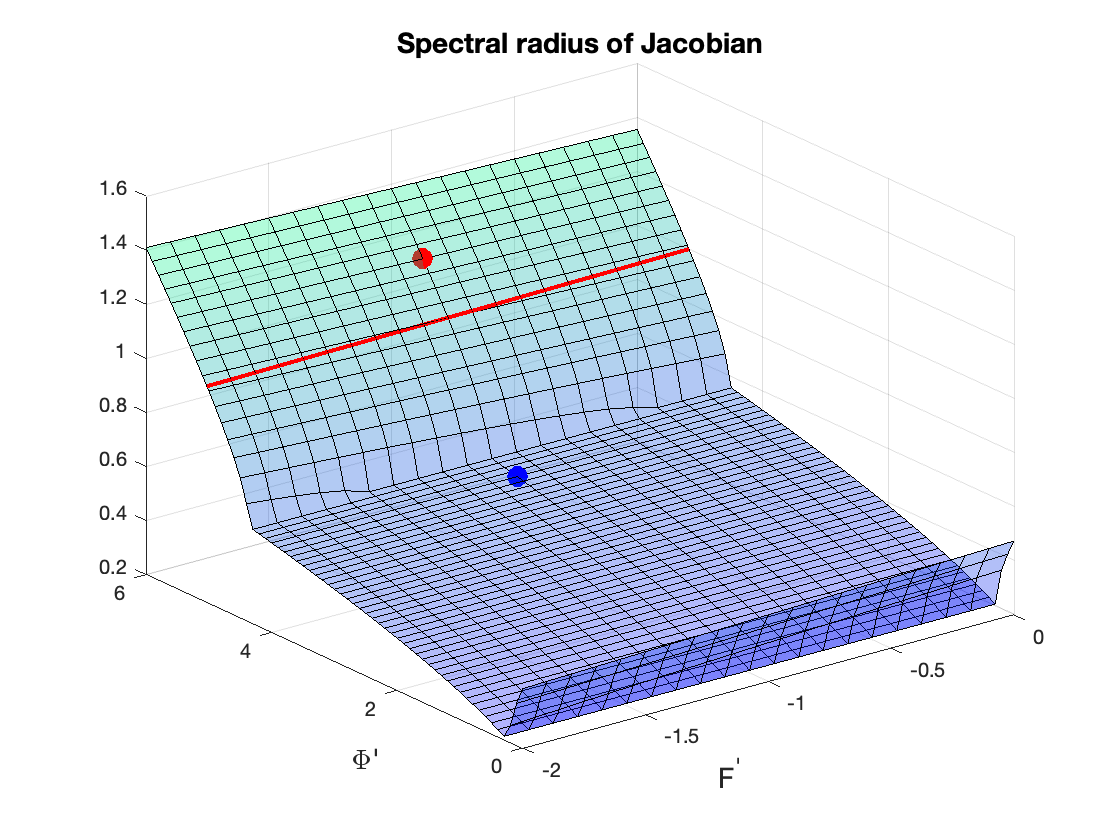}
\caption{Spectral radius of $Q^\prime(X)$ and condition \eqref{eq.stability-slopes} as function of $\Phi^\prime(y_0)$ and $F^\prime(y_0)$. Red line  -- stability border \eqref{eq.stability-slopes}. Blue dot -- stable 1-cycle with $F^\prime(y_0)=-1$, $\Phi^\prime(y_0)=4$. Red dot -- unstable 1-cycle with  $F^\prime(y_0)=-1$, $\Phi^\prime(y_0)=5.5$. }\label{fig:convergence_F_Phi}
\end{figure}

To illustrate the theoretical results of Section~\ref{sec:stability}, consider the pharmacokinetic-pharmodynamic  model of the muscle relaxant {\it atracurium} used under general  closed-loop  anesthesia. The model originates from \cite{SWM12} and has been used for investigating the performance of the IGO as a feedback dosing algorithm in \cite{MPZ24}, \cite{MPZ24_ECC}, \cite{MPZ24_MED}. For the population mean values of the patient-specific parameters, the linear part of the model is of third order (see \eqref{eq:1} with the state matrix
\[
A=\begin{bmatrix}
     -0.0374         &0         &0\\
    0.0374   &-0.1496         &0\\
         0    &0.0560   &-0.3740
\end{bmatrix}.
\]
The fixed point $X^\top=\begin{bmatrix}  269.5974 &84.5819 & 13.6249\end{bmatrix}$ corresponds to the 1-cycle with the parameters $\lambda=300$, $T=20$ see \cite{MPZ24_ECC} for details on the design procedure.
Then,
\[
J=\begin{bmatrix}
    0.4733\\
    0.1410\\
    0.0221
\end{bmatrix}, \quad D=\begin{bmatrix}
    -10.0829\\
   -2.5705\\
   -0.3633
\end{bmatrix},
\]
and inequality \eqref{eq.stability-slopes} gives the stability condition for the 1-cycle
\begin{equation}\label{eq:numerical_stabilty}
    0.0454\cdot F^\prime(y_0) -0.8550\cdot \Phi^\prime(y_0)>-1.
\end{equation}
In Fig.~\ref{fig:convergence_F_Phi}, one can see  the numerically calculated spectral radius of $Q^\prime(X)$ with the stability border of criterion \eqref{eq.stability-slopes} depicted.  The values of $F^\prime(y_0)$ and $\Phi^\prime(y_0)$  under the stability border yield stable 1-cycles that correspond to the fixed point $X$. Notice that orbital instability of a 1-cycle (i.e. unstable fixed point $X$)  means that the IGO possesses another  attractor, a stable one, whereto solutions converge when the initial condition does not coincide with the fixed point.

\begin{figure}[t]
\centering
\includegraphics[width=0.75\linewidth]{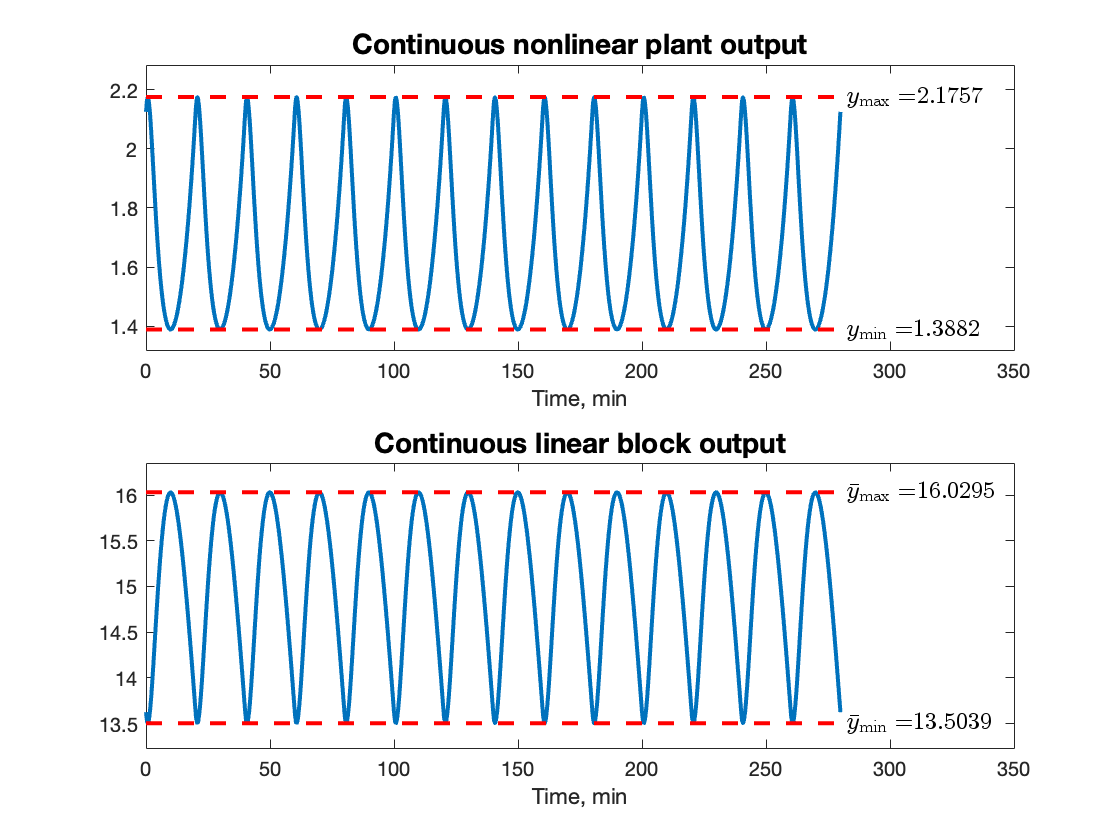}
\caption{ The  1-cycle with $\lambda=300, T=20$ stabilized by the modulation function slopes $F^\prime(y_0)=-1$, $\Phi^\prime(y_0)=4$. The initial condition on the continuous block is $x(0)=X$. Output corridor values for the linear and nonlinear output are marked.
}\label{fig:1-cycle}
\end{figure}

\begin{figure}[t]
\centering
\includegraphics[width=0.75\linewidth]{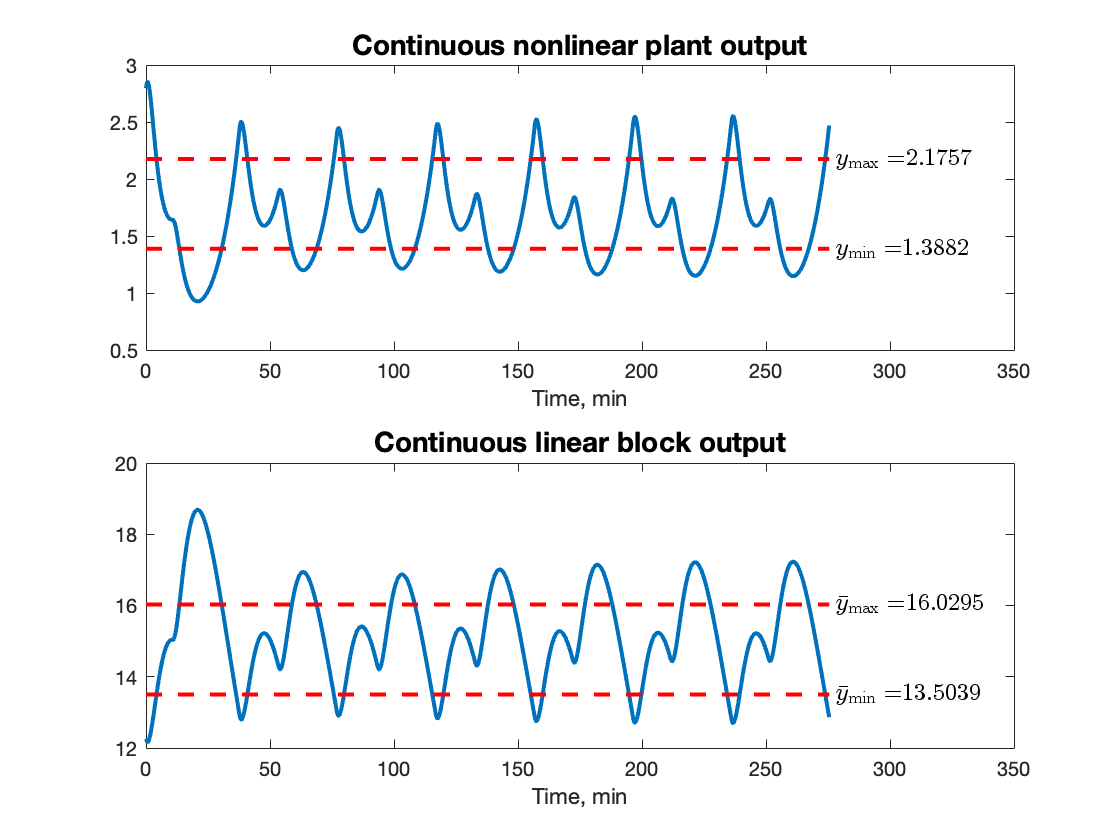}
\caption{ The  2-cycle with $\lambda=300, T=20$ stabilized by the modulation function slopes $F^\prime(y_0)=-1$, $\Phi^\prime(y_0)=5.5$. The initial condition on the continuous block is $x(0)\ne X$. Output corridor values for the linear and nonlinear output in the designed 1-cycle are marked.
}\label{fig:2-cycle}
\end{figure}

The 1-cycle with the desired parameters corresponding to the fixed point $X$ is shown in Fig.~\ref{fig:1-cycle}. When the  slopes of the modulation functions are set to  $F^\prime(y_0)=-1$, $\Phi^\prime(y_0)=5.5$ (the red dot in Fig.~\ref{fig:convergence_F_Phi}), the 1-cycle becomes orbitally unstable. Being initiated outside of the fixed point $X$, the closed-loop solution of the IGO exhibits a stable 2-cycle instead, Fig.~\ref{fig:2-cycle}. A 2-cycle is characterized by two firings of the pulse-modulated feedback over the least solution period. The occurrence of a 2-cycle  is expected as the main bifurcation mechanism in the IGO is frequency doubling \cite{ZCM12b}. This example reveals a remarkable property of the designed impulsive controller. Even though the desired solution (the 1-cycle) has become unstable, the closed-loop system output is kept in vicinity of the desired output corridor. This is in sharp contrast with instability of a equilibrium in a linear time-invariant system.

Notice that the slopes of the frequency and amplitude modulation functions control the convergence to the 1-cycle under perturbation, whereas the parameters of the stationary solution (defined by $X$) remain the same. Multiple numerical studies of the IGO's dynamics (see e.g. \cite{ZCM12b}) indicate that a stable 1-cycle attracts almost all feasible (positive) solutions.

\begin{figure}[t]
\centering
\includegraphics[width=0.75\linewidth]{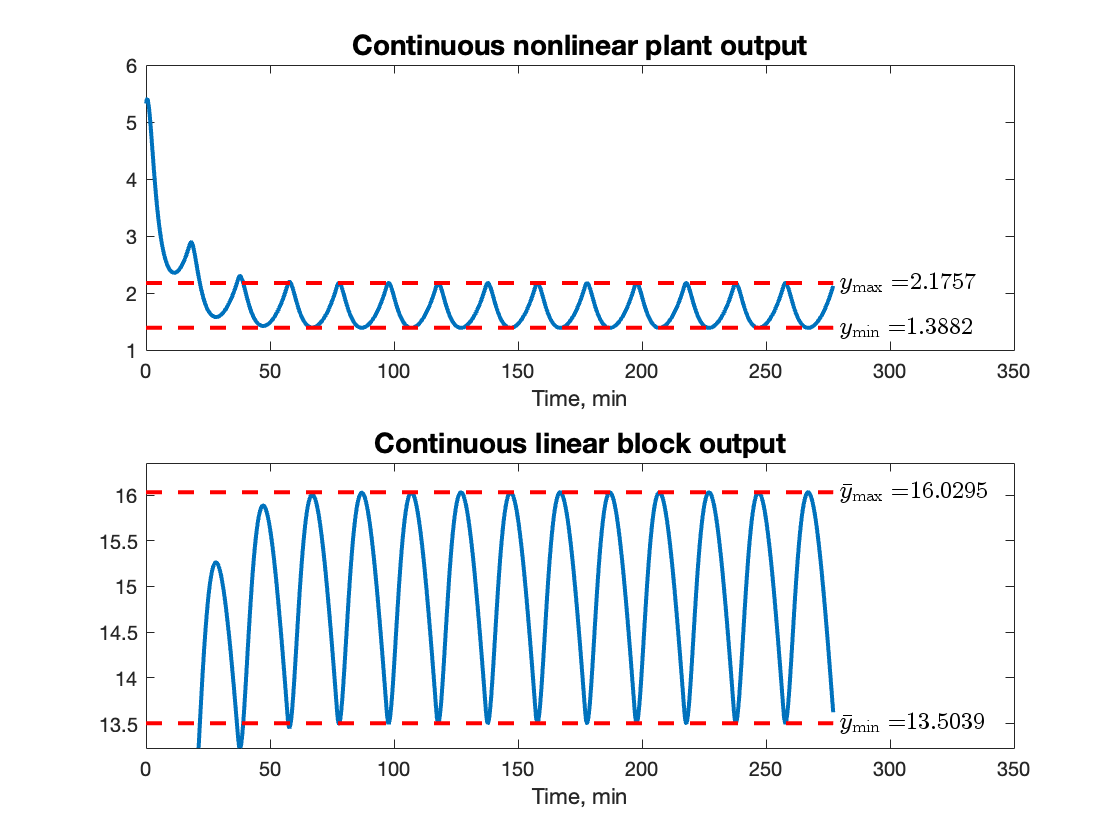}
\caption{ The  1-cycle with $\lambda=300, T=20$ stabilized by the modulation function slopes $F^\prime(y_0)=-0.1$, $\Phi^\prime(y_0)=0.29$. The initial condition on the continuous block is far from the fixed point $X$. The Jacobian eigenvalues are $ \sigma(Q^\prime(X))=\{0.2348, 0.1814, 0.0003\}$. Output corridor values for the linear and nonlinear output are marked.
}\label{fig:1-cycle-case_1}
\end{figure}

\begin{figure}[t]
\centering
\includegraphics[width=0.75\linewidth]{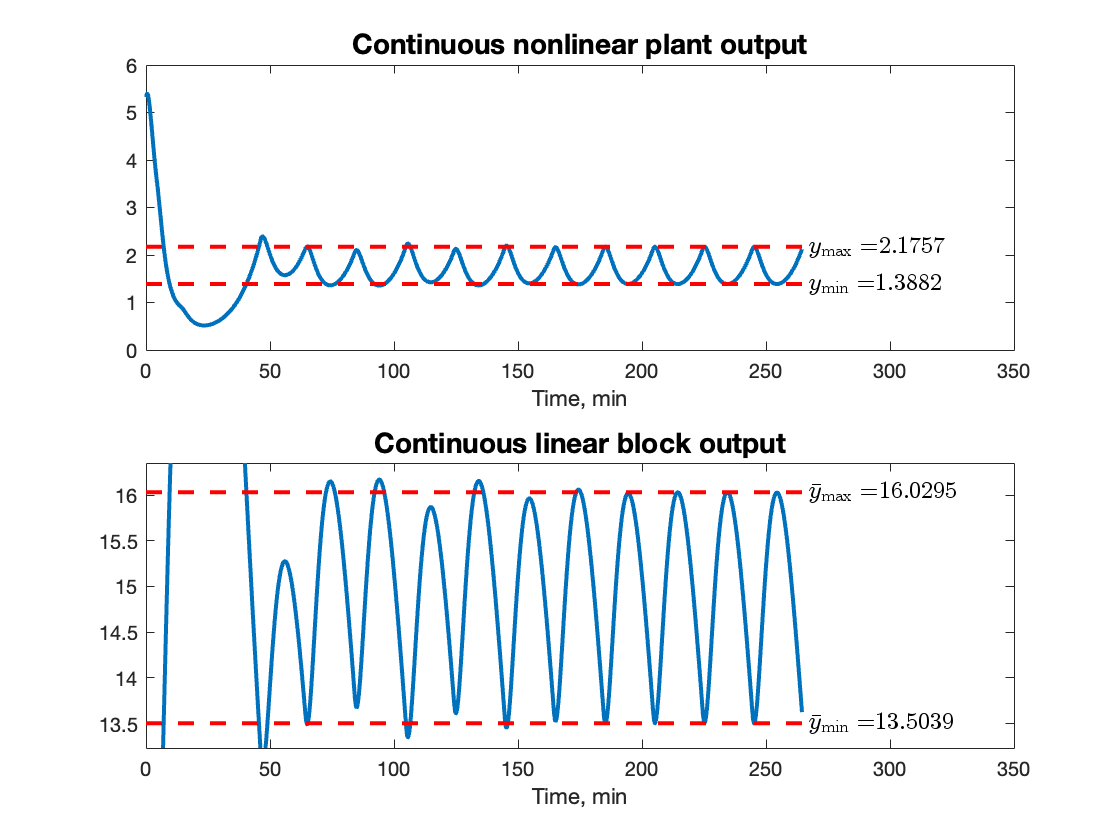}
\caption{ The  1-cycle with $\lambda=300, T=20$ stabilized by the modulation function slopes $F^\prime(y_0)=-1$, $\Phi^\prime(y_0)=4$. The initial condition on the continuous block is far from the fixed point $X$. The Jacobian eigenvalues are $ \sigma(Q^\prime(X))=\{ -0.4757 \pm 0.2343i,
  -0.0000\}$. The spectral radius is $\rho(Q^\prime(X))=0.5302$. Output corridor values for the linear and nonlinear output are marked.
}\label{fig:1-cycle-case_6}
\end{figure}

The eigenvalues of the Jacobian $Q^\prime(X)$, also known as multipliers, define the character of the transient process to the stationary periodic solution. To illustrate this, consider an IGO design with $F^\prime(y_0)=-0.1$, $\Phi^\prime(y_0)=0.29$. In Fig.~\ref{fig:1-cycle-case_1}, the transient response of the closed-loop system to the stationary solution is simulated. The eigenvalues of the Jacobian are all real, positive, and less than one ($ \sigma(Q^\prime(X))=\{0.2348, 0.1814, 0.0003\}$), which results in a monotonic (in envelope) convergence to the stationary solution.
Fig.~\ref{fig:1-cycle-case_6} illustrates a case of the Jacobian spectrum including a complex pair, namely $ \sigma(Q^\prime(X))=\{ -0.4757 \pm 0.2343i,
  -0.0000\}$. The transient response shows a strong overshoot that is undesirable in many applications, e.g. drug dosing and process control.

\section{CONCLUSIONS}

Stability of the 1-cycle in the Impulsive Goodwin's Oscillator  (IGO) is examined. A linear inequality specifying the stability domain of the stationary solution in terms of the slopes of the frequency and amplitude modulation functions is derived. The result is shown to be instrumental in optimizing the convergence rate and dynamical character of perturbed solutions to the 1-cycle under stability guarantee. A formal solution of the optimization problem is saved for future works. The IGO gives rise to a class of simple feedback controllers that implement administration of discrete doses to a continuous plant according to a desired schedule or maintain the plant output in a given closed interval of values. The presence of the pulse-modulated feedback allows the IGO to manipulate both the doses and their timing to achieve the control goal.

\section*{ACKNOWLEDGMENT}
A. Medvedev was partially supported by the Swedish Research Council under grant 2019-04451.


\section*{APPENDIX: Proof of Theorem~\ref{thm.main}}

Throughout this section, all the assumptions of Theorem~\ref{thm.main} are supposed to hold, in particular, $0<a_1<a_2<a_3$.
Consider the matrix
\begin{equation}\label{eq:Q_general}
        \begin{gathered}
        \mathcal{Q}(T,\xi,\eta)= \e^{AT}+ \left( \xi J+\eta D\right)C,\\
        D\triangleq A(\e^{-AT}-I)^{-1}B,\,J=\e^{AT}B,
        \end{gathered}
    \end{equation}
where $A,B,C$ are matrices from~\eqref{eq:1}.

We first prove a technical lemma.

\begin{lemma}\label{lem.non-complex-conj}
Suppose that $T>0$, $\xi\leq 0$, and $\eta\geq 0$. Then matrix $\mathcal Q\triangleq\mathcal{Q}(T,\xi,\eta)$ has the following spectral properties:
\begin{enumerate}
\item  $\mathcal Q$ has no eigenvalues on the interval $(\e^{-a_{1}T},\infty)$;
\item there exists a real eigenvalue $z_1\in [\e^{-a_3T},\e^{-a_{1}T}]$;
\item the product of two remaining eigenvalues $z_2z_3$ does not exceed $e^{-(a_1+a_2)T}<1$.
\item $\mathcal{Q}$ is not Schur stable if and only if $z_2,z_3$ are real and $\min(z_2,z_3)\leq -1$;
\end{enumerate}
\end{lemma}
\proof
The characteristic polynomial $\chi(z)\triangleq \det(zI-\mathcal{Q})$, thanks to the Schur complement formula, is written as
\begin{equation}\label{eq.aux1}
\begin{aligned}
\chi(z)&=\det(zI-\e^{AT})w(z),\quad\text{where}\\
w(z)&\triangleq\frac{\chi(z)}{\det(zI-\e^{AT})}=\\
&=1-C(zI-\e^{AT})^{-1}(\xi J+\eta D).
\end{aligned}
\end{equation}
Notice that $w(z)\geq 1$ whenever $z$ is real and $z>\e^{-a_{1}T}$, because $1-w(z)$ can be decomposed into the series
\[
C(zI-e^{AT})^{-1}\left(\xi J+\eta D\right)=z^{-1}\sum_{k=0}^{\infty}z^{-k}Ce^{kAT}\left(\xi J+\eta D\right),
\]
whose coefficients are non-positive matrices, because $J>0$ and $D<0$~\cite[Proposition~3]{MPZh23}.
Furthermore, from the triangular structure of $A$, the characteristic polynomial
\begin{equation}\label{eq.aux2}
\det(zI-\e^{AT})=(z-\e^{-a_1T})(z-\e^{-a_2T})(z-\e^{-a_3T})
\end{equation}
is positive for all $z>\e^{-a_{1}T}$, which entails that $\chi(z)>0$. This completes the proof of statement 1).

To prove statement 2), a more subtle argument is needed that requires us to compute the residual of rational function $w(z)$ at $z=\e^{-a_3T_0}$.
To this end, consider the diagonalization of matrix $A$. It can be checked that
\[
\begin{gathered}
S^{-1}AS=
\begin{bmatrix}
-a_1 & 0 & 0\\
0 & -a_2 & 0\\
0 & 0 & -a_3
\end{bmatrix},\;\;\text{where}\\
S=
\begin{bmatrix}
1 & 0 & 0\\
\tfrac{g_1}{a_2-a_1} & 1 & 0\\[3pt]
\tfrac{g_1g_2}{(a_2-a_1)(a_3-a_1)} & \tfrac{g_2}{a_3-a_2} & 1
\end{bmatrix}\;\;\text{and}\\
S^{-1}=
\begin{bmatrix}
1 & 0 & 0\\
-\tfrac{g_1}{a_2-a_1} & 1 & 0\\[3pt]
-\tfrac{g_1g_2}{(a_3-a_2)(a_3-a_1)} & -\tfrac{g_2}{a_3-a_2} & 1
\end{bmatrix}.
\end{gathered}
\]
Denote for brevity $\bar B\triangleq S^{-1}B$, $\bar C\triangleq CS$, that is,
\[
\bar B=\begin{bmatrix}
1\\
-\tfrac{g_1}{a_2-a_1}\\[3pt]
-\tfrac{g_1g_2}{(a_3-a_2)(a_3-a_1)}
\end{bmatrix},\;
\bar C=
\begin{bmatrix}
\tfrac{g_1g_2}{(a_2-a_1)(a_3-a_1)} \\ \tfrac{g_2}{a_3-a_2} \\ 1
\end{bmatrix}^{\top}
\]
Considering the function
\[
\rho_z(s)=(zI-\e^{Ts})^{-1}\left(\xi \e^{Ts}+\eta(\e^{-Ts}-1)^{-1}s\right),
\]
one obtains the relation
\[
\begin{aligned}
\rho_z(A)&=S\rho_z(S^{-1}AS)S^{-1}=\\&=S
\begin{bmatrix}
\rho_z(-a_1) & 0 & 0\\
0 & \rho_z(-a_2) & 0\\
0 & 0 & \rho_z(-a_3)
\end{bmatrix}
S^{-1}
\end{aligned}
\]
Hence, the rational function $w$ from~\eqref{eq.aux1} can be written as
\[
w(z)=1-\sum_{i=1}^3\bar c_i\bar b_i\rho_z(-a_i).
\]
It can now be noticed that the residual of $w$ at $z=\e^{-a_3T}$ is non-positive.
Indeed, $\bar c_3=1>0,\bar b_3<0$, and hence
\[
\begin{aligned}
\lim_{z\to \e^{-a_{3}T}}&\big(z-\e^{-a_{3}T}\big)w(z)=\\&=-\bar c_3\bar b_3\lim_{z\to \e^{-a_{3}T}}\big(z-\e^{-a_{3}T}\big)\rho_z(-a_3)=\\
&=-\bar c_3\bar b_3\left(\xi \e^{-a_3T}-\eta a_3(\e^{a_3T}-1)^{-1}\right)\leq 0.
\end{aligned}
\]
Here we used the fact functions $\rho_{z}(\e^{-a_1T})$, $\rho_{z}(\e^{-a_2T})$ are analytic in $z$ in a vicinity of $\e^{-a_3T}$, and also $\xi\leq 0,\eta\geq 0$.
On the other hand, recalling the definition of $w(z)$ and~\eqref{eq.aux2},
the latter residual can be computed as
\[
\begin{aligned}
\lim_{z\to \e^{-a_{3}T}}&\big(z-\e^{-a_{3}T}\big)w(z)=\\&=
\frac{\det\big(\e^{-a_3T}I-\mathcal{Q}\big)}{(\e^{-a_3T}-\e^{-a_2T})(\e^{-a_3T}-\e^{-a_1T})},
\end{aligned}
\]
entailing that
\[
\det(zI-\mathcal{Q})|_{z=\e^{-a_{3}T}}\leq 0.
\]
At the same time, it has been already proven that
\[
\det(zI-\mathcal{Q})|_{z=\e^{-a_{1}T}}\geq 0,
\]
which implies that second statement.

To prove the remaining statements, it suffices to notice that
\[
\begin{aligned}
\det\mathcal{Q}&=\det\e^{AT}\left(1+C\e^{-AT}\left( \xi J+\eta D\right)\right)=\\
&=\det\e^{AT}\left(1+\xi CB+\eta CA(I-\e^{AT})^{-1}B\right)
\end{aligned}
\]
Obviously, $CB=0$. It can be shown that the function
\[
z\mapsto \psi(z)\triangleq\frac{z}{1-e^{z}}
\]
is concave on the interval $z\in (-\infty,0)$.

Using the Opitz formula (see, e.g., Step~2 in the proof of~\cite[Lemma~11]{PRM24}), one obtains that
\[
\begin{aligned}
CA(I-e^{AT})^{-1}B=T^{-1}C\psi(TA)B=\\=T^{-1}\psi[-a_1T,-a_2T,-a_3T].
\end{aligned}
\]
The generalized mean-value theorem~\cite[Lemma~10]{PRM24} entails now the existence of $\zeta\in(-a_3T,-a_1T)$ such tha $\psi[-a_1T,-a_2T,-a_3T]=\psi''(\zeta)/2$.
Thanks to the concavity of $\tau$, one thus has $CA(I-e^{AT})^{-1}B\leq 0$, whence
\[
z_1z_2z_3=\det\mathcal{Q}\leq \det\e^{AT}=\e^{-(a_1+a_2+a_3)T}.
\]
This implies statement 3) in virtue of $z_1\geq \e^{-a_3T}$.

To prove statement 4), it suffices to notice that a pair of complex-conjugate eigenvalues $z_2=z_3^*$ should have the modulus
$|z_2|=|z_3|\leq e^{-(a_1+a_2)T/2}<1$. Hence, if $\mathcal{Q}$ has one real and two complex-conjugate eigenvalues, it is automatically Schur stable. The only reason for being unstable is
thus the existence of a \emph{real} eigenvalue whose modulus is not less than $1$. In view of statement 1), $\mathcal{Q}$ cannot have eigenvalue at $1$. Hence, one of $z_2,z_3$ does not exceed $-1$
(in which case the remaining eigenvalue is, obviously, also real).
\QED
\begin{corollary}\label{cor.Schur-stab}
Let the assumption of Lemma~\ref{lem.non-complex-conj} apply. Then, the following three statements are equivalent:
\begin{enumerate}
\item Matrix $\mathcal Q\triangleq\mathcal{Q}(T,\xi,\eta)$ is Schur stable;
\item The inequality holds as follows
\begin{equation}\label{eq.stability-det0}
\chi(-1)=\det(-I-\mathcal{Q})<0,
\end{equation}
\item $\xi,\eta$ obey the inequality
\begin{equation}\label{eq.stability-slopes0}
C(I+\e^{TA})^{-1}\left(\xi J+\eta D\right)>-1.
\end{equation}
\end{enumerate}
\end{corollary}
\proof
To prove that conditions~\eqref{eq.stability-det0} and~\eqref{eq.stability-slopes0} are equivalent, it suffices to substitute $z=-1$ into~\eqref{eq.aux1}
and notice that
\[
\det(-I-\e^{AT})=-(1+\e^{-a_1T})(1+\e^{-a_2T})(1+\e^{-a_3T)})<0.
\]
Hence,~\eqref{eq.stability-slopes0} holds (equivalently, $w(-1)>0$) if and only if $\chi(1)<0$, i.e., statements 2) and 3) are equivalent.

Obviously, 1) implies 2), because $\chi(z)\to-\infty$ when $z$ is real and $z\to-\infty$. If one has $\chi(-1)\geq 0$, then matrix $\mathcal{Q}$ has an eigenvalue on $(-\infty,-1]$ and is thus not Schur stable.

To prove that 2) (and 3)) implies 1), consider a one-parameter family of matrices
$\mathcal{Q}_{\varepsilon}\triangleq\mathcal{Q}(T,\varepsilon\xi,\varepsilon\eta)$ and the corresponding characteristic polynomials $\chi_{\varepsilon}(z)\triangleq\det(zI-\mathcal{Q}_{\varepsilon})$.
Notice that if~\eqref{eq.stability-slopes0} holds, then it remains valid by replacing $\xi,\eta$ by $\varepsilon\xi,\varepsilon\eta$, where $\varepsilon\in[0,1]$.
Hence, $\chi_{\varepsilon}(-1)<0$ for all $\varepsilon\in[0,1]$.
Obviously, $\chi_0$ is a Schur polynomial and the coefficients of $\chi_{\varepsilon}$ continuously depend on $\varepsilon$. Hence, either $\chi_{\varepsilon}$ is Schur for all $\varepsilon\in[0,1]$, or
there exists $\varepsilon_0$ such that $\chi_{\varepsilon_0}$ has a root on the unit circle $\mathbb{S}=\{z\in\mathbb{C}:|z|=1\}$.
The second alternative is, however, impossible: Lemma~\ref{lem.non-complex-conj}, applied to $\mathcal{Q}(T,\varepsilon_0\xi,\varepsilon_0\eta)$,
states that the only possible eigenvalue on $\mathbb{S}$ is $z=-1$, whereas $\chi_{\varepsilon_0}(-1)<0$.
\QED

\subsection*{The proof of Theorem~\ref{thm.main}}

The proof is straightforward from Corollary~\ref{cor.Schur-stab}. Applying the latter Corollary~\ref{cor.Schur-stab} to $T=\Phi(y_0)$, $\xi=F'(y_0)$, $\eta=\lambda\Phi'(y_0)$, one easily checks that $Q'(X)=\mathcal{Q}(T,\xi,\eta)$, where the matrix-valued function $\mathcal{Q}$ is defined in~\eqref{eq:Q_general}.\QED

\bibliographystyle{IEEEtran}
\bibliography{observer,refs}

\end{document}